\documentclass[twocolumn]{article}
\usepackage{graphicx} 
\usepackage{abstract} 
\usepackage{natbib}
\usepackage{amsmath}
\usepackage{amssymb}
\usepackage{algpseudocode}
\usepackage{algorithm}
\usepackage{color}

\title{Approximate Dynamic Programming for a Remanufacture-to-Order System}
\author{Amirreza Pashapour\thanks{Graduate School of Science and Engineering, Koç University, İstanbul, Türkiye (e-mail: apashapour20@ku.edu.tr)}, Fatemeh Zare Bidaki\thanks{Graduate School of Business, Koç University, İstanbul, Türkiye (e-mail: fbidaki19@ku.edu.tr)}}
\date{}
\begin{document}

\maketitle

\begin{abstract}
    Remanufacturing is pivotal in transitioning to more sustainable economies. While industry evidence highlights its vast market potential and economic and environmental benefits, remanufacturing remains underexplored in theoretical research. This study revisits a state-of-the-art remanufacture-to-order (RtO) system and develops an alternative approach by developing an approximate dynamic programming (ADP) algorithm to solve larger RtO instances. The proposed methodology yields results consistent with existing state-dependent, non-congestive policies. Key findings include the optimality of fulfilling demand whenever remanufacturable cores are available and prioritizing cores of the highest quality level. This work highlights the potential of ADP in addressing problems in remanufacturing domains.
\end{abstract}

\textbf{Keywords:} remanufacturing, inspection, Markov decision process, dynamic programming, approximate dynamic programming, optimal control

\section{Introduction}

Remanufacturing plays a pivotal role in driving economies and supply chains toward a more sustainable and circular model. Remanufactured products typically offer a cost advantage of 40–60\% compared to their new parts \citep{munot2015insight}, achieve an environmental impact reduction of 80\% in energy and water consumption, and decrease waste generation by 70\% \citep{Perkins2019}. Remanufacturing firms rebuild worn-out products, typically known as cores, by replacing worn components or reprocessing parts. This process stands out as a key strategy in restoring end-of-life (EOL) products to meet their original quality standards regarding specifications and warranty \citep{matsumoto2016trends}. Examples of such cores include automotive components, and air conditioning compressors and motors, for which the practice from the industry has shown that there is a huge market \citep{atasu2008remanufacturing}. 

Within remanufacturing systems, tactical challenges primarily arise in the product acquisition and demand fulfillment phases \citep{souza2013closed}. Acquisition quantity is a critical managerial decision for remanufacturers, as it directly impacts unit remanufacturing costs \citep{galbreth2010optimal,mutha2016managing,abbey2019improving}. While determining the optimal acquisition quantity is straightforward when the conditions of used products are known with certainty, remanufacturers often face inherent ambiguity during acquisition \citep{nadar2023optimal}. The quality of returned products varies significantly, influenced by factors such as user behavior and operational environments \citep{li2016core}. This variability can result in batches of cores with predominantly low quality. Although remanufacturers may possess prior knowledge of the quality distribution, the exact condition of products is typically unknown at the point of acquisition. Inspection and testing are essential to mitigate this uncertainty and enable firms to classify used-items into quality grades or cores, each associated with distinct holding and remanufacturing costs \citep{ferguson2009value,fleischmann2010product}. Following this classification, remanufacturers must also make optimal decisions regarding fulfilling an arrived demand and, if so, determine the appropriate core to utilize.

Our study is motivated by the observation that, despite its significant potential to benefit both the environment and the economy, remanufacturing remains underexplored in theoretical research \citep{matsumoto2016trends}. A state-of-the-art study by \cite{nadar2023optimal} examines the used-item acquisition and disposition problem within a single-product RtO system involving unknown quality conditions for acquired cores. In this context, the actual quality type of a used-item is revealed only after acquisition and inspection. Any unmet demand is lost. Demand for the remanufactured product can be fulfilled using any available core; however, cores with higher quality incur lower remanufacturing costs.

To address this problem, \cite{nadar2023optimal} formulate the problem as a Markov Decision Process (MDP) and incorporate the uncertainty associated with used-item conditions at the time of acquisition. Used-item procurement lead times follow an exponential distribution, demand for remanufactured products arrives according to a Poisson process, and the system state is defined by the inventory levels of each core type. Their work identifies a state-dependent non-congestive optimal policy that, for each system state, specifies whether to procure a used-item, whether to fulfill an incoming demand, and, if so, which core to remanufacture to satisfy it. They further established key structural properties for the optimal cost function from a congestion control perspective. The authors tested their results on instances having two quality types. Building on this foundation, our goal in the present work is to exercise an ADP algorithm on the RtO system studied by \cite{nadar2023optimal}. Upon designing an ADP algorithm, our work revisits and extends their analysis to larger problem instances while validating the proposed model, results, and findings through the structural properties established in the original study. 

The remainder of this paper is structured as follows. Section 2 provides an overview of the related literature. Section 3 revisits the RtO system, drawing on the problem description and formulation presented by \cite{nadar2023optimal}. In Section 4, we develop an ADP framework specifically tailored for the RtO system. Section 5 presents the design of baseline experimental instances of varying sizes, along with their corresponding results and findings, which aligns with the findings analyzed in \cite{nadar2023optimal}. Finally, Section 6 concludes the paper.

\section{Literature review}

The remanufacturing literature is well-studied, with works addressing diverse aspects of product acquisition and disposition. Studies dealing with product acquisition can be classified into three main directions: (i) passive acceptance of all returns and subsequent classification (e.g. \cite{gong2013optimal,inderfurth2013advanced,jia2016addressing}; (ii) pricing-driven strategies to influence the quality and quantity of returns (e.g. \cite{karakayali2007analysis,zhou2011optimal,cai2014optimal}; and (iii) procurement without pricing flexibility, where firms must optimize acquisition and disposition decisions under fixed market conditions (e.g. \cite{teunter2011optimal,panagiotidou2013optimal,vercraene2014coordination,li2016core,mutha2016managing,gayon2017optimal,liu2022remanufacturing,nadar2023optimal}. Our study examines a remanufacturing system similar to that evaluated by \cite{nadar2023optimal}, which falls within the third category and focuses on acquisition decisions in a fixed market environment with condition uncertainty for a RtO system. We approach the problem from a different methodological perspective by developing an ADP framework.

The analysis of (re)manufacturing systems using dynamic programming (DP) is also well-established in the literature. In this area, most studies focus on product-level management, typically for a single component. \cite{fleischmann2003optimal} pioneered in investigating a single-component inventory system with returns and showed the optimality of a stationary $(s, S)$ order policy. \cite{teunter2008multi} addressed a scheduling problem of two production sources—manufacturing and remanufacturing—and found that dedicated production is more cost-efficient than integrating both operations on a single line. \cite{zhou2011optimalOR} explored product acquisition, pricing, and inventory management decisions in remanufacturing, while \cite{zhou2011optimalMSOM} extended this to multiple core types returned from customers. \cite{kim2013joint} studied a single-component system integrating production, remanufacturing, and disposal decisions and showed that the optimal policy follows a state-dependent base-stock structure. \cite{vercraene2014coordination} analyzed scenarios where remanufactured products are perceived as identical to new ones, where the optimal policy is given by state-dependent base-stock thresholds. \cite{cai2014optimal} focused on core quality and identified a state-dependent base-stock policy as optimal. \cite{gayon2017optimal} investigated the disposal of returns, either immediately upon arrival or after they become serviceable. \cite{panagiotidou2017joint} further considered product quality in returned items. Recently, \cite{nadar2023optimal} developed an MDP model for the used-item acquisition and disposition problem within a single-product RtO system involving multiple unknown quality conditions for acquired cores. \cite{dekker2013reverse} and \cite{khan2021state} have offered thorough, systematic reviews of remanufacturing systems and their potential operational difficulties. Our contributions include addressing the option of directly discarding cores with inferior qualities within the RtO model and developing and exercising an ADP algorithm tailored for a RtO system.

\section{Problem definition: An MDP model}

The problem we overview in this section is driven based on the recent study of \cite{nadar2023optimal}. Consider a used-item acquisition problem faced by a remanufacturer dealing with uncertain conditions of worn-out items of a single product. At the time of acquisition, the true physical condition of each core is unknown. However, the remanufacturer has prior knowledge of the distribution of these conditions, consistent across all acquired cores. Once a core is acquired, it undergoes disassembly for part inspection, during which its actual quality, categorized into one of $\mathcal{K}=\{1,\ldots,K\}$ types, is determined. Let the core types be ranked from highest quality $(i=1)$ to lowest $(i=K)$. The remanufacturer maintains an inventory for each core type $1$ to $K$. The unit cost of remanufacturing of these cores $r_i$ varies with core quality such that $r_1 \leq r_2 \leq \ldots \leq r_K$. As an extension to \cite{nadar2023optimal}'s model, we incorporate a $K+1$ type for uneconomical cores that can be detected with probability $\overline{p}$ and must be discarded. Prior to an inspection, the remanufacturer knows that each core will fall into type-i with a probability of $p_i$ for $i\in\mathcal{K}\cup \{K+1\}$, a common assumption in the remanufacturing literature \citep{vedantam2021revenue}. Let the unit selling price of the remanufactured product be exogenous to the system and be greater than $r_i$ for $i \in \mathcal{K}$, thereby remanufacturing any core be profitable.  

The acquisition and demand arrivals are characterized as follows. The remanufacturer can acquire and inspect at most one core at a time which will incur an acquisition cost of $c_a$. Consistent with typical inventory models, we assume the acquisition lead times are independent and exponentially distributed with mean $1/\mu$. The available cores are remanufactured on demand, and no inventory for remanufactured products is maintained. RtO systems are well-documented in the literature \citep{teunter2011optimal,jia2016addressing,li2016core} and are widely adopted by various industries, including U.S. Naval Aviation Depots, large-scale mailing equipment manufacturers, and cell phone remanufacturers. Demand for the remanufactured product follows a Poisson process with a rate of $\lambda$. If no cores are available in inventory, demand must be rejected. Additionally, demand may be declined even when cores are in stock. Each unmet demand incurs a unit lost sales cost, $c_l$, where $c_l \geq r_K$.

The system state at any time $t$ is the vector $\mathbf{x}_t=(x_{1t},\ldots,x_{Kt})$, where the nonnegative integer $x_{it} \in \mathbb{N}^0$ denotes the number of type $i$ cores at time $t$. We consider a capacity bound $b$ on the number of cores in the inventory, i.e., $\sum_{i \in \mathcal{K}} x_{it} \leq b$. Each core $i$ held in stock incurs a unit holding cost $h_i > 0$ per unit time where $h_1 \geq h_2 \geq \ldots \geq h_K$. Let $h(\mathbf{x}_t)=\sum_{i \in \mathcal{K}} h_i x_{it}$ be the total inventory holding cost rate once in state $\mathbf{x}_t$. Let $t_0 = 0$, and $t_j$ be the time of the $j^{\text{th}}$ state transition. Since all inter-event times follow an exponential distribution, the system exhibits the memoryless property. As a result, the problem is modeled as an MDP, the optimality equation remains consistent for all $t \in [t_j, t_{j+1}]$, and decision-making can be limited to the specific moments when the system's state changes.

A control policy $\pi$ specifies, for each state $\mathbf{x}$, the action $u^{\pi}(\mathbf{x}) = \{\tau, \eta\}$, where $\tau \in \{0, 1\}$ and $\eta \in \{0, 1, \dots, K\}$. The $\tau = 1$ indicates that the acquisition process is initiated; otherwise, $\tau = 0$. Similarly, if $\eta \neq 0$, demand for the remanufactured product is satisfied using a type-$\eta$ core provided that at least one type-$\eta$ core is available in inventory; if $\eta = 0$, demand is rejected. The set of admissible actions for state $\mathbf{x}$ is denoted by $\mathcal{U}(\mathbf{x})$. Let $V(.)$ be a real-valued function $\mathbb{R}^K \rightarrow \mathbb{R}$ that denotes the expected discounted cost over an infinite planning horizon and $\alpha \in (0, 1)$ be the discount parameter. 

The objective is to identify a policy $\pi^*$ that minimizes the expected discounted costs. Without loss of generality, the time scale is set to $\lambda + \mu = \alpha$ to embed the effect of discounting within the optimality equation. With this scale, we formulate the optimality equation that holds for the optimal cost function $V^* = V^{\pi^*}$ as in Equation \ref{eq:optimality}.
\begin{equation} \label{eq:optimality}
    V(\mathbf{x}) = \mu T_1 V(\mathbf{x}) + \lambda T_2 V(\mathbf{x}) + h(\mathbf{x})
\end{equation}

The operator $T_1$ in Eq. \ref{eq:operatorT1} specifies whether to initiate the acquisition of a used-item given $\mathbf{x}$. The operator $T_2$ in Eq. \ref{eq:operatorT2} decides whether to fulfill an arrived demand as well as which core type to remanufacture if the demand is to be satisfied. In Eq. \eqref{eq:operators}, the $e_i$ is the $i^{\text{th}}$ unit vector in dimension $K$ for $i \in \mathcal{K}$. Also, $\overline{p}$ is the probability that the inspected item is inferior and should be discarded. Note that $\overline{p} + \sum_{i \in \mathcal{K}} p_i = 1$. 
\begin{small}
\begin{subequations} \label{eq:operators}
\begin{eqnarray}
    T_1 V(\mathbf{x}) = \min \left\{V(\mathbf{x}), c_a + \sum_{i \in \mathcal{K}} p_i V(\mathbf{x}+e_i) + \overline{p} V(\mathbf{x}) \right\} \label{eq:operatorT1} \\ 
    T_2 V(\mathbf{x}) = \min \left\{ V(\mathbf{x}) + c_l , \min_{i \in \mathcal{K}:\; x_i \geq 1} V(\mathbf{x}-e_i) + r_i\right\} \;\;\; \label{eq:operatorT2}
\end{eqnarray}
\end{subequations}
\end{small}
Given $b$ as the inventory capacity of all cores and $K$ the core qualities, the size of the state space is $|X|=\binom{b+K}{K}=(b+K)!/(b!K!)$. Also, the size of the action space is $|\mathcal{U}|=2K+2$. We observe that although the action space is tractable, the state space grows combinatorially.

\cite{nadar2023optimal} explored structural properties for the optimal cost function \eqref{eq:optimality}, established a state-dependent noncongestive-acquisition policy, and designed a numerical experiment testbed considering two core types for the optimal policy characterization. Diverging from \cite{nadar2023optimal}'s solution approach, in the present work, we adopt a different methodology to tackle larger RtO instances with more core types. For this purpose, we develop an ADP framework as detailed in Section 4.

\section{Approximate Dynamic Programming}

We propose an ADP algorithm to obtain high-quality core acquisition and demand satisfaction policies for RtO. In this regard, we establish parametric linear approximations of the value functions by defining a set of state-dependent basis functions $\phi(\mathbf{x})$ in Eq. \eqref{eq:basis} (See \cite{ccil2009effects} for some operators in queuing systems). The goal is to accurately calibrate the $\theta$-vector of size $K+1$ that represents the weights of each basis function, 
\begin{equation} \label{eq:basis}
 \overline{V}(\mathbf{x}) \approx \theta^{\top} \phi(\mathbf{x}) = \theta_0 g(\mathbf{x}) + \sum_{i \in \mathcal{K}} \theta_i \frac{x_i}{b},
\end{equation}
where,
\begin{equation} \label{eq:gx}
 g(\mathbf{x}) = \left(\frac{\lambda}{\lambda + \mu}\right)^{\mathbf{x}^{\top}\mathbf{1}}
\end{equation}
is a $g:\mathbb{R}^{K} \rightarrow \mathbb{R}$ function. We utilize $g(\mathbf{x})$ to effectively address the critical challenge of depleting the core inventory required to meet demand. This function is a decreasing convex function which makes it well-suited to capture the trade-offs inherent in inventory control. Specifically, when inventory levels are critically low, $g(\mathbf{x})$ reflects the heightened risk of demand loss. The terms $x_i/b$ in Eq. \eqref{eq:basis} serve as the rest of our basis functions, which incorporate key cost factors such as holding and remanufacturing costs. We note that all basis functions are scaled within the range of 0 to 1.

Our proposed ADP leverages an Approximate Policy Iteration (API) algorithmic framework. Based on exact Policy Iteration (PI), the API starts with a base policy, that is, an initial $\theta$-vector, and iteratively generates a sequence of improved $\theta$-vectors and accordingly approximate value functions. This process operates within two nested loops: an outer loop for policy improvement (PI) and an inner loop for policy evaluation (PE) \citep{bertsekas2011approximate}. The loops are explained as follows.

Given a $\theta$-vector by the outer loop, the inner loop evaluates the performance of the current policy. During each inner loop iteration $z \in \{ 1, \ldots, Z \}$, a feasible state $\mathbf{x}_z$ is first randomly sampled, and its basis function evaluation $\phi(\mathbf{x}_z)$ is recorded. Subsequently, using Eq. \eqref{eq:operators}, we acquire the best course of action, i.e. $u^{* \theta}(\mathbf{x})=\{\tau^*, \eta^*\}$, for the $T_1$ and $T_2$ operators based on a $\epsilon$-greedy rule. In this rule, with a probability of $\epsilon$, we explore the admissible action space by selecting an action at random, whereas with a probability of $1 - \epsilon$, we exploit the optimal course of action by applying the arg-min operator to $T_1$ and $T_2$.

Now, it is time to simulate the exogenous information $w$ in the system. In particular, this involves the occurrence of (i) $w=\rho_a$, a core acquisition (with rate $\mu$) or (ii) $w=\rho_d$, a demand fulfillment (with rate $\lambda$) event. Given the Markovian property of these events, their probability of occurrence is $\mu/(\mu + \lambda)$ and $\lambda/(\mu + \lambda)$, respectively. Let $t_w$ be the duration of time elapsed until the occurrence of either event. According to the revealed information, we compute and record the immediate cost contribution of being in state $\mathbf{x}$ and taking action $u^{* \theta}(\mathbf{x})$ using Eq. \eqref{eq:immediatecost}. Finally, we simulate the system for one event forward and record the basis function evaluation $\phi(\mathbf{x'}_z)$.
\begin{subequations} \label{eq:immediatecost}
\begin{equation} \label{eq:immediatecost1}
C^{\rho_{a}}\left(\mathbf{x}, u^{* \theta}(\mathbf{x})\right) = t_w h(\mathbf{x}) + 
\begin{cases}
    0 & \text{if } \tau = 0 \\
    c_a & \text{if } \tau = 1 \\
\end{cases}
\end{equation}
\begin{equation} \label{eq:immediatecost2}
C^{\rho_{d}}\left(\mathbf{x}, u^{* \theta}(\mathbf{x})\right) = t_w h(\mathbf{x}) +
\begin{cases}
    c_l & \text{if } \eta = 0 \\
    r_{\eta} & \text{if } \eta \geq 1 \\
\end{cases}
\end{equation}
\end{subequations} 

Upon the termination of the inner loop, $Z$ observations are recorded into matrices $\Omega$ and $\Omega'$, as our basis function evaluations, both with dimensions $Z \times (K+1)$, along with a column vector of immediate costs $\hat{C}$ of size $Z \times 1$, as outlined in Eq. \eqref{eq:matrices}.
\begin{equation} \label{eq:matrices}
\Omega = 
\begin{bmatrix}
\phi(\mathbf{x}_1) \\ \phi(\mathbf{x}_2) \\ \vdots \\ \phi(\mathbf{x}_Z)
\end{bmatrix}, \;
\Omega' = 
\begin{bmatrix}
\phi(\mathbf{x}'_1) \\ \phi(\mathbf{x}'_2) \\ \vdots \\ \phi(\mathbf{x}'_Z)
\end{bmatrix}, \;
\hat{C} = 
\begin{bmatrix}
C\left(u^{* \theta}(\mathbf{x}_1)\right) \\ C\left(u^{* \theta}(\mathbf{x}_2)\right) \\ \vdots \\ C\left(u^{* \theta}(\mathbf{x}_Z)\right)
\end{bmatrix}
\end{equation}

The outer loop leverages this recorded data to calibrate $\theta$. The process begins by computing $\hat{\theta}$ as an estimation of $\theta$, through the Least-Squares Temporal Difference (LSTD) learning framework. Originally proposed by \citet{bradtke1996linear}, LSTD is a computationally efficient method to estimate the tunable vector$\theta$, particularly when the value functions are approximated using a linear combination of the basis functions. The LSTD performs a least squares regression, minimizing the sum of squares of the $Z$ temporal differences $\theta^{\top}\phi(\mathbf{x}_z) + C\left(\mathbf{x}, u^{* \theta}(\mathbf{x})\right) - \alpha \theta^{\top}\phi(\mathbf{x}'_z)$.

Given the matrices $\Omega$ and $\Omega'$ and the vector $\hat{C}$, we utilize Eq. \eqref{eq:updateTheta}, also known as the ``Instrumental variables method" \citep{soderstrom2002instrumental}, to attain $\hat{\theta}$. In this equation, we incorporate a $(K+1)\times (K+1)$ diagonal regularization matrix $\beta I$. The diagonal entries are controlled by the regularization parameter $\beta \geq 0$. Regularization prevents singularity during matrix inversion and diminishes generalization errors.
\begin{equation} \label{eq:updateTheta}
 \hat{\theta}  =  {\left[\Omega^{\top}(\Omega - \alpha \Omega') + \beta I \right]}^{-1} \; {(\Omega)}^{\top} \hat{C}
\end{equation}

Ultimately, the outer loop updates the $\theta$-vector using a harmonic step-size rule as in Eq. \eqref{eq:thetaupdate}. On its right-hand side, $\theta$ denotes an estimate based on all previous outer loop iterations, and $\hat{\theta}$ is the current estimate of the outer loop iteration $n=1,\ldots, N$. The coefficient $\gamma_n = 1/n^\delta$ is a polynomial smoothing step-size coefficient, where $\delta \in (0.5,1]$. The weight of individual estimates reduces, and antecedent estimates are prioritized as $n$ rises.
\begin{equation} \label{eq:thetaupdate}
 \theta \gets (1 - \gamma_n) \theta + \gamma_n \hat{\theta}
\end{equation}
The steps of the proposed ADP are provided in Algorithm \ref{alg:mainADP}. The algorithm accurately estimates the $\theta$-vector as a high-quality policy vector for all states. The tunable algorithmic parameters include $N$ (the number of PI iterations), $Z$ (the number of PE iterations), $\beta$ (regularization term), and $\delta$ (polynomial step-size parameter).

{\centering
\begin{algorithm}[H] \footnotesize
 \caption{API Algorithm}
 \label{alg:mainADP}
 \begin{algorithmic}[1]
  \State{Initialize the $\theta$-vector}
  \For{$n=1$ to $N$}: \hfill \Comment{Policy Improvement loop}
   \For{$z=1$ to $Z$}: \hfill \Comment{Policy Evaluation loop}
    \State{Pick a feasible random state $\mathbf{x}$}
    \State{Record the basis function evaluation $\phi(\mathbf{x})$}
    \State{Obtain the best decision $u^{* \theta}(\mathbf{x})$ via a $\epsilon$-greedy rule}
    \State{Simulate $w$: a core acquisition or demand arrival event}
    \State{Record the associated cost contribution}
    \State{Given $\mathbf{x}$, $w$ and $u^{* \theta}(\mathbf{x})$, simulate to a new state $\mathbf{x}'$}
    \State{Record the basis function evaluation $\phi(\mathbf{x}')$}
    \EndFor
   \State{\textbf{Apply LSTD:} Compute $\hat{\theta}$}
   \State{Update $\theta$}
   \EndFor
  \State{\textbf{return} $\theta$ and $\overline{V}(.|\theta)$}
 \end{algorithmic}
\end{algorithm}
\par}
In the following section, we conduct numerical experiments to evaluate the performance of the proposed ADP method on RtO instances of varying sizes.

\section{Numerical results}

\subsection{Design of RtO instances} \label{sec:design}

\citet{nadar2023optimal} employed a standard value iteration algorithm on a testbed of 243 instances with $K=2$. In this section, we design a new testbed of size 12 ($K \in \{2,3,4,5\}$ and $\lambda \in \{0.25,0.5,0.75\}$) with a focus on larger state variables. The baseline experimental design is presented in Table \ref{tab:param}.
\begin{table}[ht]
\footnotesize
    \centering
    \begin{tabular}{cccc}
    \multicolumn{2}{c}{RtO factors} & \multicolumn{2}{c}{Algorithmic factors} \\ \hline
    Parameter & Value & Parameter & Value \\ \hline
    $K$ & \{2, 3, 4, 5\} & $N$ & 10 \\
    $h_i$ & $K - i + 1$ & $Z$ & 1000 \\
    $r_i$ & $10i$ & $\theta$ & $[1,h_1,\ldots,h_K]$ \\
    $b$ & 20 & $\beta$ & 10 \\
    $\lambda$ & \{0.25, 0.5, 0.75\} & $\delta$ & 0.5 \\
    $\mu$ & 0.99 $- \lambda$ & $\epsilon$ & 0.05 \\
    $c_a$ & 5 & $\alpha$ & 0.99 \\
    $c_l$ & 100 & & \\
    $p_i$ & $1/(K+1)$ & & \\
    \end{tabular}
    \caption{Baseline parameter values for our computational experiments}
    \label{tab:param}
\end{table}

\subsection{Results and findings}
We implemented the proposed ADP algorithm in Python. Each of the 12 instances is solved 10 times to account for the inherent randomness in the process. Consequently, the results and findings presented in this section are based on the average $\theta$ vector obtained across these repetitions.

\begin{figure}[ht]
\centering
\includegraphics[width=0.8\linewidth]{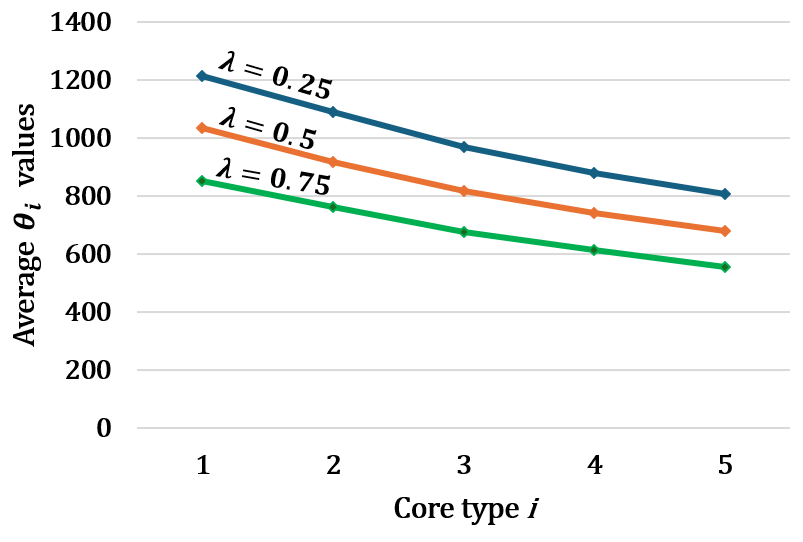}
\caption{\centering Average values of $\theta_i$ coefficients for core types $i \in \mathcal{K}$}
\label{fig:theta}
\end{figure}

\begin{itemize}
 \item $\theta_0 > 0$ and $\theta_i > 0 \; \forall i \in \mathcal{K}$: Consistent with the optimal policy structure outlined in \cite{nadar2023optimal}, we observed that it is optimal to fulfill demand whenever cores are available for remanufacturing. This arises from the fact that the $\theta$ vector is revealed to be elementwise non-negative in the outputs.
 
 \item $\theta_0 \nearrow \lambda$: We observed that $\theta_0$ roughly converged to values of 129, 206, and 415 for $\lambda \in \{0.25, 0.5, 0.75\}$. This indicates that as the demand arrival rate increases, the system incurs a higher risk of lost sales in the case of low core inventory.

 \item $\theta_i \searrow \lambda \; \forall i \in \mathcal{K}$: Unlike $\theta_o$, we observed that $\theta_i$ coefficients for each core $i \in \mathcal{K}$ revealed a decreasing trend with respect to $\lambda$, as seen in Figure \ref{fig:theta}. This is also a reasonable observation since as the demand rate rises, core inventories, in the long run, incur fewer holding costs. 
 
 \item $\theta_i \searrow i \; \forall \lambda$: We observed that $\theta_1 \geq \theta_2 \geq \ldots \geq \theta_K$ (see Figure \ref{fig:theta}) which aligns with the optimal policy structure discussed in the literature. It indicates that it is always optimal to prioritize remanufacturing cores of the highest quality level. This is because higher-quality cores impose greater costs on the system over time, making it more beneficial to process them first.

\end{itemize}

The above results illustrate the effective performance of the designed algorithm. Once the $\theta$ vector has converged, we obtain the approximate value functions. By incorporating these approximate values into Eq. \eqref{eq:optimality}, we can quickly determine the best course of action via an $\arg \min (.)$ function. In this context, we aim to compute thresholds that serve as order-up-to levels. 

\begin{figure}[ht]
\centering
\includegraphics[width=0.8\linewidth]{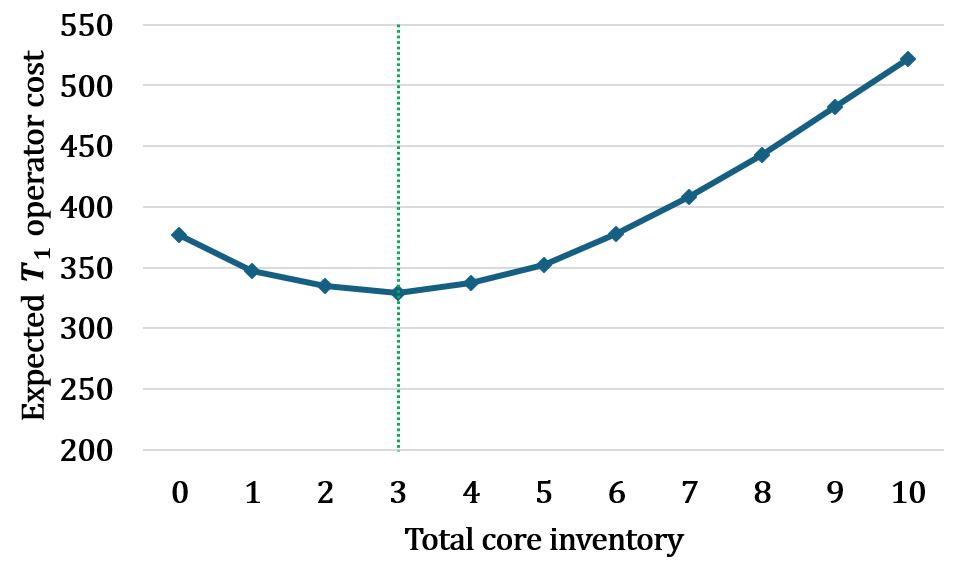}
\caption{\centering Expected $T_1$ operator costs associated with different total core inventory levels}
\label{fig:t1}
\end{figure}

\begin{itemize}
\item \textit{When to acquire:} We seek to understand if we are in state $\mathbf{x}$ and a used-item arrives, whether it is optimal to acquire it or not. For this, we analyze Eq. \eqref{eq:operatorT1} under varying state $\mathbf{x}$ for the case of $K=5$ and $\lambda = 0.75$. Particularly, we first enumerate all possible states whose components sum to $\overline{x}$, compute their corresponding $T_1$ operator costs, and take the average of these values. The results for total core inventories $\overline{x}$ are presented in Figure \ref{fig:t1}. Given the baseline setting in Section \ref{sec:design}, as long as the total core inventory is less than or equal to 3, it is optimal to acquire an arriving used-item; otherwise, it is optimal to reject it. This result is consistent with the proposed state-dependent noncongestive-acquisition policy proposed by \cite{nadar2023optimal}. Additionally, we note that due to the relationship $\lambda + \mu = \alpha$ in the proposed system, the proportion of demand or used-item arrivals is relatively low. If this proportion were higher, the threshold could increase to accommodate larger inventory levels.
\end{itemize}

\section{Conclusion}

Remanufacturing plays a pivotal role in advancing our economy and supply chains toward more sustainable and circular models. Remanufactured products typically offer a cost advantage of 40–60\% compared to new parts, reduce environmental impact by up to 80\% in energy and water consumption, and decrease waste generation by 70\%. This process involves rebuilding worn-out products, or cores, by replacing damaged components or reprocessing parts. However, remanufacturing systems face significant tactical challenges, particularly in the acquisition of used products and the fulfillment of demand. Among these, the quantity of cores acquired is a critical managerial decision, as it directly affects remanufacturing costs.  

Despite its dual potential to benefit the environment and the economy, remanufacturing remains underexplored in theoretical research. A recent state-of-the-art study by \cite{nadar2023optimal} addresses the used-item acquisition and disposition problem in a single-product RtO system with multiple unknown quality conditions for acquired cores. The authors formulate this problem as an MDP by incorporating uncertainties associated with core conditions at the time of acquisition. They propose a state-dependent, non-congestive optimal policy that specifies whether to procure a core, whether to fulfill demand and, if so, which core to remanufacture to meet the demand.  

Building on this foundation, our work employs an ADP algorithm to extend the analysis of the RtO system studied by \cite{nadar2023optimal}. We adopt an ADP algorithm for larger RtO instances. Our analysis confirms the optimality of prioritizing high-quality cores for remanufacturing and consistently fulfilling demand whenever cores are available. Moreover, we observe the persistence of the order-up-to structure in ADP results for larger state sizes.  

In conclusion, we hope our study not only validates the structural properties established in prior work but also paves the way for broader applications of ADP in remanufacturing system studies.

{\small

}

\end{document}